\documentclass[11pt]{article}

\newcommand{\makefig}[3]{
	\begin{figure}[htbp]
        \refstepcounter{figure}
	\label{#2}
        \begin{center}
		~#3~\\
		\medskip
                {\sf Figure \thefigure.  #1}
        \end{center}
	\medskip
	\end{figure}
}

\edef\atcatcode{\the\catcode`\@}
\catcode`@=11 

\newif\ifoldlabels \oldlabelsfalse
\expandafter\ifx\csname @ifundefined\endcsname\relax 
\long\def\@ifundefined#1#2#3{\expandafter\ifx\csname
  #1\endcsname\relax#2\else#3\fi}
\def\@iden#1{#1}
\def\@warning#1{\message{Warning: #1.}}
\fi

\@ifundefined{ps@init}{}{\endinput}

\@ifundefined{hyperactivelabels}{\let\hyperactivelabels\relax}{}

\newif\if@topspecials 
\input driver.chg

\newcount\psfc@a
\newcount\psfc@b
\newcount\psfc@c
\newcount\psfc@d
\newcount\psfc@e
\newcount\psfc@f
\newcount\psfc@g
\newcount\psfc@h
\newcount\psfc@i
\newcount\psfc@j
\newcount\psscale@count
\newtoks\psfigtoks@
\newdimen\psfig@dimen
\def\@nnil{\@nil}
\def\@empty{}
\def\@psdonoop#1\@@#2#3{}
\def\@psdo#1:=#2\do#3{\edef\@psdotmp{#2}\ifx\@psdotmp\@empty \else
    \expandafter\@psdoloop#2,\@nil,\@nil\@@#1{#3}\fi}
\def\@psdoloop#1,#2,#3\@@#4#5{\def#4{#1}\ifx #4\@nnil \else
       #5\def#4{#2}\ifx #4\@nnil \else#5\@ipsdoloop #3\@@#4{#5}\fi\fi}
\def\@ipsdoloop#1,#2\@@#3#4{\def#3{#1}\ifx #3\@nnil 
       \let\@nextwhile=\@psdonoop \else
      #4\relax\let\@nextwhile=\@ipsdoloop\fi\@nextwhile#2\@@#3{#4}}
\def\@tpsdo#1:=#2\do#3{\xdef\@psdotmp{#2}\ifx\@psdotmp\@empty \else
    \@tpsdoloop#2\@nil\@nil\@@#1{#3}\fi}
\def\@tpsdoloop#1#2\@@#3#4{\def#3{#1}\ifx #3\@nnil 
       \let\@nextwhile=\@psdonoop \else
      #4\relax\let\@nextwhile=\@tpsdoloop\fi\@nextwhile#2\@@#3{#4}}
\def\psdraft{%
  \def\@psdraft{0}
}
\def\psfull{%
  \def\@psdraft{100}
}
\psfull
\newif\if@prologfile
\newif\if@postlogfile
\newif\if@noisy
\def\pssilent{%
  \@noisyfalse
}
\def\psnoisy{%
  \@noisytrue
}
\pssilent
\newif\if@bbllx
\newif\if@bblly
\newif\if@bburx
\newif\if@bbury
\newif\if@height
\newif\if@width
\newif\if@rheight
\newif\if@rwidth
\newif\if@clip
\newif\if@proportional
\newif\if@scale
\newif\if@verbose
\def\@p@@sclip#1{\@cliptrue}
\def\@p@@sproportional#1{\@proportionaltrue}
\def\@p@@sfile#1{%
  \def\@p@sfile{#1.ps}%
  \def\@labfile{#1.lab}%
}
\def\@p@@sfigure#1{%
  \def\@p@sfile{#1.ps}%
  \def\@labfile{#1.lab}%
}
\def\@pDimenToSpNumber #1#2{%
  \psfig@dimen = #2\relax
  \edef#1{\number\psfig@dimen}%
}
\def\@p@@sbbllx#1{%
  \@bbllxtrue
  \@pDimenToSpNumber{\@p@sbbllx}{#1}%
}
\def\@p@@sbblly#1{%
  \@bbllytrue
  \@pDimenToSpNumber{\@p@sbblly}{#1}%
}
\def\@p@@sbburx#1{%
  \@bburxtrue
  \@pDimenToSpNumber{\@p@sbburx}{#1}%
}
\def\@p@@sbbury#1{%
  \@bburytrue
  \@pDimenToSpNumber{\@p@sbbury}{#1}%
}
\def\@p@@sheight#1{%
  \@heighttrue
  \@pDimenToSpNumber{\@p@sheight}{#1}%
}
\def\@p@@swidth#1{%
  \@widthtrue
  \@pDimenToSpNumber{\@p@swidth}{#1}%
}
\def\@p@@srheight#1{%
  \@rheighttrue
  \@pDimenToSpNumber{\@p@srheight}{#1}%
}
\def\@p@@srwidth#1{%
  \@rwidthtrue
  \@pDimenToSpNumber{\@p@srwidth}{#1}%
}
\def\@p@@ssilent#1{%
  \@verbosefalse
}
\def\@p@@sscale #1{%
  \def\@p@scale{#1}%
  \@scaletrue
}

\def\@p@@sprolog#1{\@prologfiletrue\def\@prologfileval{#1}}
\def\@p@@spostlog#1{\@postlogfiletrue\def\@postlogfileval{#1}}
\def\@cs@name#1{\csname #1\endcsname}
\def\@setparms#1=#2,{\def\@tempa{#2}\ifx\@tempa\@empty
  \@warning{bad syntax (missing = or extra comma) in argument of \string\psfig}%
  \else\@@setparms#1=#2,\fi}
\def\@@setparms#1=#2=,{\@cs@name{@p@@s#1}{#2}}
\def\ps@init@parms{%
  \@bbllxfalse \@bbllyfalse
  \@bburxfalse \@bburyfalse
  \@heightfalse \@widthfalse
  \@rheightfalse \@rwidthfalse
  \@scalefalse
  \def\@p@sbbllx{}\def\@p@sbblly{}%
  \def\@p@sbburx{}\def\@p@sbbury{}%
  \def\@p@sheight{}\def\@p@swidth{}%
  \def\@p@srheight{}\def\@p@srwidth{}%
  \def\@p@sfile{}%
  \def\@labfile{}%
  \def\@p@scost{10}%
  \def\@sc{}%
  \@prologfiletrue
  \@postlogfilefalse
  \@clipfalse
  \@proportionalfalse
  \if@noisy
    \@verbosetrue
  \else
    \@verbosefalse
  \fi
}
\def\parse@ps@parms#1{%
   \@psdo\@psfiga:=#1\do
     {\expandafter\@setparms\@psfiga=,}%
}

\newif\ifno@file
\newif\ifno@bb
\newif\ifnot@eof
\newif\if@bbmatch
\newif\if@crematch
\newif\ifmathematica
\newif\ifillustrator
\newread\ps@stream
\newread\lab@stream
\newcount\@linecount
\newcount\maxheaderlines
\def\scan@header{%
  \openin\ps@stream=\@p@sfile
  \openin\lab@stream=\@labfile
  \ifeof\ps@stream
    \relax 
    \@warning{cannot open \@p@sfile}
    \no@filetrue
  \else
    \not@eoftrue
    \ifno@bb \@bbmatchfalse \else \@bbmatchtrue \fi
    \@crematchfalse
    \catcode`\%=12
    \catcode`\:=12 
    \@linecount=\maxheaderlines
    \loop
      \read\ps@stream to \line@in
      \global\psfigtoks@=\expandafter{\line@in}
      \ifeof\ps@stream \not@eoffalse \fi
      \if@bbmatch \else \@bbtest{\psfigtoks@} \fi
      \if@crematch \else \@cretest{\psfigtoks@} \fi
      \if@crematch \if@bbmatch \not@eoffalse \fi \fi
      \advance\@linecount-1
      \ifnum\@linecount=0 \not@eoffalse \fi
    \ifnot@eof \repeat
  \fi
  \catcode`\%=14
}  

{ 
  \catcode`\%=12
  \catcode`\:=12
  \gdef\@bbtest#1{\expandafter\@bb@\the#1
  \global\long\def\@bb@#1
    \else\@bbmatchtrue\expandafter\bb@cull\the\psfigtoks@\fi}

  \gdef\@cretest#1{\expandafter\@cre@\the#1
  \global\long\def\@cre@#1
    \else\@crematchtrue\@mathtest{\psfigtoks@}\@illtest{\psfigtoks@}\fi}
}

\def\@mathtest#1{\expandafter\@math@\the#1Mathematica\@mathtest\@math@}
\long\def\@math@#1Mathematica#2#3\@math@{\ifx\@mathtest#2
    \else\mathematicatrue \fi}

\def\@illtest#1{\expandafter\@ill@\the#1Illustrator\@illtest\@ill@}
\long\def\@ill@#1Illustrator#2#3\@ill@{\ifx\@illtest#2
    \else\illustratortrue \fi}

\def\bb@cull#1:{\expandafter\bb@@cull\@iden}
\long\def\bb@@cull#1 #2 #3 #4 {%
  \@pDimenToSpNumber{\@p@sbbllx}{#1bp}%
  \@pDimenToSpNumber{\@p@sbblly}{#2bp}%
  \@pDimenToSpNumber{\@p@sbburx}{#3bp}%
  \@pDimenToSpNumber{\@p@sbbury}{#4bp}%
  \no@bbfalse
}

\def\compute@bb{%
  \no@bbfalse
  \if@bbllx \else \no@bbtrue \fi
  \if@bblly \else \no@bbtrue \fi
  \if@bburx \else \no@bbtrue \fi
  \if@bbury \else \no@bbtrue \fi
  \scan@header 
  \ifno@file
  \else
    \ifno@bb
      \@warning{no bounding box found in \@p@sfile}
      \no@filetrue
    \else
      \psfc@c=\@p@sbburx
      \psfc@b=\@p@sbbury
      \advance\psfc@c by -\@p@sbbllx
      \advance\psfc@b by -\@p@sbblly
      \edef\@bbw{\number\psfc@c}
      \edef\@bbh{\number\psfc@b}
    \fi
  \fi
}

\def\in@hundreds #1#2#3{%
  \psfc@g=#2
  \psfc@d=#3
  \divide\psfc@d 10
  \psfc@a=\psfc@g  
  \divide\psfc@a by \psfc@d
  \psfc@f=\psfc@a
  \multiply\psfc@f by \psfc@d
  \advance\psfc@g by -\psfc@f
  \multiply\psfc@g by 10
  \psfc@f=\psfc@g  
  \divide\psfc@f by \psfc@d
  \psfc@j=\psfc@f
  \multiply\psfc@j by \psfc@d
  \advance\psfc@g by -\psfc@j
  \multiply\psfc@g by 10
  \psfc@j=\psfc@g  
  \divide\psfc@j by \psfc@d
  \psfc@h=#1\psfc@i=0
  \psfc@e=\psfc@h
  \multiply\psfc@e by \psfc@a
  \advance\psfc@i by \psfc@e
  \psfc@e=\psfc@h
  \divide\psfc@e by 10
  \multiply\psfc@e by \psfc@f
  \advance\psfc@i by \psfc@e
  \psfc@e=\psfc@h
  \divide\psfc@e by 100
  \multiply\psfc@e by \psfc@j
  \advance\psfc@i by \psfc@e
  \divide\psfc@i 10
  \edef\@result{\number\psfc@i}
}
\def\@ScaleInHundreds #1{%
  \in@hundreds{#1}{\@p@scale}{100}%
  \edef#1{\@result}%
}
\def\compute@wfromh{%
  \in@hundreds{\@p@sheight}{\@bbw}{\@bbh}%
  \edef\@p@swidth{\@result}%
}
\def\compute@hfromw{%
  \in@hundreds{\@p@swidth}{\@bbh}{\@bbw}%
  \edef\@p@sheight{\@result}%
}
\def\compute@minhw{%
  \in@hundreds{\@p@swidth}{\@bbh}{\@bbw}%
  \ifnum\@p@sheight<\@result
    \in@hundreds{\@p@sheight}{\@bbw}{\@bbh}%
    \edef\@p@swidth{\@result}%
  \else
    \edef\@p@sheight{\@result}%
  \fi
}
\def\compute@handw{%
  \if@height 
    \if@width
      \if@proportional
        \compute@minhw
      \fi
    \else
      \compute@wfromh
    \fi
  \else 
    \if@width
      \compute@hfromw
    \else
      \edef\@p@sheight{\@bbh}
      \edef\@p@swidth{\@bbw}
    \fi
  \fi
}
\def\compute@resv{%
  \if@rheight \else \edef\@p@srheight{\@p@sheight} \fi
  \if@rwidth \else \edef\@p@srwidth{\@p@swidth} \fi
}
\def\partest{\par}
\def\psfig#1{%
  \vbox {%
    \offinterlineskip
    \ps@init@parms
    \parse@ps@parms{#1}%
    \compute@bb
                \ifno@file\vbox{\hbox{{\tt\@p@sfile} not found}}%
                \else
    \compute@handw
    \compute@resv
    \if@scale
      \@ScaleInHundreds{\@p@swidth}%
      \@ScaleInHundreds{\@p@sheight}%
      \@ScaleInHundreds{\@p@srwidth}%
      \@ScaleInHundreds{\@p@srheight}%
    \fi
    \ifnum\@p@scost<\@psdraft
      \if@topspecials\do@specials\fi
      \vbox to \@p@srheight sp{%
        \hbox to \@p@srwidth sp{}%
        \vss
      }%
      \if@topspecials\else\do@specials\fi
      \ifeof\lab@stream
      \else{%
        \hyperactivelabels
        \not@eoftrue
        \loop
          \read\lab@stream to \line@in \ifx\line@in\partest\else\line@in\fi
          \ifeof\lab@stream \not@eoffalse \fi
        \ifnot@eof \repeat
      }\fi
    \else
      \vbox to \@p@srheight sp{%
        \hbox to \@p@srwidth sp{%
          \if@verbose
            \@p@sfile
          \fi
        }%
        \vss
      }%
    \fi
  \fi
  }%
}

\def\do@specials{\ps@init\ps@begin
  \if@clip 
    \ps@clip \fi
  \if@prologfile \ps@prolog \fi \ps@include \if@postlogfile \ps@postlog \fi
  \ps@end}

\newbox\label@box
\newdimen\x@lab \newdimen\y@lab
\newdimen\x@aux \newdimen\y@aux
\newdimen\hair\hair=3pt
\def\setlabel#1#2#3#4#5{%
  \setbox\label@box\hbox{$#1$}%
  \x@lab.5\wd\label@box \x@lab#4\x@lab
  \y@lab.5\ht\label@box\advance\y@lab.5\dp\label@box \y@lab#5\y@lab 
  \x@aux.92388\x@lab \advance\x@aux.38268\y@lab
  \y@aux-.38268\x@lab \advance\y@aux.92388\y@lab
  \ifdim\x@aux>0pt
    \ifdim\y@aux>0pt
      \ifdim\x@aux>\y@aux\advance\x@lab.7071\hair\advance\y@lab.7071\hair
      \else\advance\y@lab\hair\fi
    \else
      \ifdim\x@aux>-\y@aux\advance\x@lab\hair
      \else\advance\x@lab.7071\hair\advance\y@lab-.7071\hair\fi
    \fi
  \else
    \ifdim\y@aux>0pt
      \ifdim\x@aux>-\y@aux\advance\x@lab-.7071\hair\advance\y@lab.7071\hair
      \else\advance\x@lab-\hair\fi
    \else
      \ifdim\x@aux>\y@aux\advance\y@lab-\hair
      \else
        \ifdim\x@aux<0pt 
           \advance\x@lab-.7071\hair\advance\y@lab-.7071\hair\fi
      \fi
    \fi
  \fi
  \advance\x@lab.5\wd\label@box
  \advance\y@lab.5\ht\label@box\advance\y@lab.5\dp\label@box
  \x@aux=#2bp \ifoldlabels \else \advance\x@aux by -\@p@sbbllx sp \fi
  \y@aux=#3bp \ifoldlabels \else \advance\y@aux by -\@p@sbblly sp \fi
  \in@hundreds{\x@aux}{\@p@swidth}{\@bbw}
  \edef\@xpos{\@result}
  \in@hundreds{\y@aux}{\@p@sheight}{\@bbh}
  \edef\@ypos{\@result}
  \vbox to 0pt{%
    \vss\hbox to\@p@srwidth sp{\hskip \@xpos sp \hskip-\x@lab 
    \box\label@box\hss}\vskip \@ypos sp \vskip-\y@lab}}

\catcode`@=\atcatcode

\usepackage{amssymb}

\newfam\meuffam
\font\tenmeuf=eufm10
\font\sevenmeuf=eufm7
\font\fivemeuf=eufm5

\textfont\meuffam=\tenmeuf
\scriptfont\meuffam=\sevenmeuf
\scriptscriptfont\meuffam=\fivemeuf

\newfam\msbfam
\font\tenmsb=msbm10
\font\sevenmsb=msbm7
\font\fivemsb=msbm5

\textfont\msbfam=\tenmsb
\scriptfont\msbfam=\sevenmsb
\scriptscriptfont\msbfam=\fivemsb

\def\Bbb{\fam\msbfam\tenmsb}

\def\l{{\germ l}}

\def\C{{\Bbb C}}

\def\P{{\Bbb P}}

\def\R{{\Bbb R}}

\def\T{{\Bbb T}}

\def\Z{{\Bbb Z}}

\def\sA{{\cal A}}

\def\sC{{\cal C}}

\def\sE{{\cal E}}

\def\sG{{\cal G}}

\def\sL{{\cal L}}
\def\sM{{\cal M}}

\def\sP{{\cal P}}
\def\sQ{{\cal Q}}

\def\sT{{\cal T}}

\def\sV{{\cal V}}

\newcommand\qed{\nopagebreak[4]\begin{flushright}\rule{0.1in}{0.1in}
\end{flushright}\pagebreak[2]}

\def\degree{\mathrm{degree}}

\font\l=cmr10 at 10pt
\font\ls=cmr7
\font\lss=cmr5
\font\lsy=cmsy10
\font\lsys=cmsy7
\font\lsyss=cmsy5
\font\lmi=cmmi10
\font\lmis=cmmi7
\font\lmiss=cmmi5
\font\lex=cmex10

 at 6truept
 at 6truept

\def\mapright#1{\smash{
  \mathop{\longrightarrow}\limits^{#1}}}
\def\mapdown#1{\Big\downarrow
  \rlap{$\vcenter{\hbox{$\scriptstyle#1$}}$}}

\newcommand\cd[1]{\matrix{#1}}

\newtheorem{theorem}{Theorem}
\newtheorem{lemma}{Lemma}[section]
\newtheorem{proposition}[lemma]{Proposition}
\newtheorem{corollary}[lemma]{Corollary}
\newcommand\heading[1]{\smallskip\noindent{\bf
#1}}

\title{Boundary Manifolds of Line Arrangements}
\author{Eriko Hironaka \thanks{Research partially supported by
N.S.E.R.C. grant OGP0170260}}
\begin{document}
\maketitle
\begin{abstract}
In this paper we describe the complement of real line
arrangements in the complex plane in terms of the boundary 
three-manifold of the line arrangement.
We show that the boundary manifold of any line arrangement
is a graph manifold with Seifert fibered vertex manifolds, 
and depends only on the incidence graph of the arrangement. 
When the line arrangement is defined over the real numbers,
we show that the homotopy type of  the complement is determined 
by the incidence graph together with orderings on the edges 
emanating from each vertex.    
\end{abstract}

\section{Introduction}
Let $\sL$ be a finite union of lines in the complex plane $\C^2$.
Two line arrangements $\sL_1$ and $\sL_2$ are said to be {\it 
topologically equivalent} if there is a homeomorphism of pairs
$$
(\C^2,\sL_1)  \rightarrow (\C^2,\sL_2).
$$
The {\it incidence graph} $\Gamma_\sL$ associated to $\sL$ is the 
bipartite graph with line-vertices corresponding to the set of lines
$\sA$ in $\sL$, point-vertices corresponding to the set of points of
intersection $\sP$ on $\sL$ and edges $e(p,L)$ and $e(L,p)$ 
whenever $p \in L$.  A morphism between incidence graphs is
a morphism of graphs preserving the vertex labelings, so that
line-vertices go to line-vertices and point-vertices go to 
point-vertices.  Two line arrangements $\sL_1$ and $\sL_2$
are said to be {\it combinatorially equivalent} if there
is an isomorphism of incidence matrices $\Gamma_{\sL_1} 
\rightarrow \Gamma_{\sL_2}$. 

Our motivation is to understand,
given a line arrangement $\sL$ in the complex plane $\C^2$, 
to what extent the topology of the pair $(\C^2,\sL)$
is determined by the combinatorics of $\sL$.  It is known 
that the cohomology of the complement $E_\sL = \C^2 \setminus 
\sL$ only depends on the incidence graph $\Gamma_\sL$
\cite{O-T:Arr}, \cite{G-M:Strat}, while 
the homotopy type and fundamental group of $E_\sL$ depend on more
information \cite{Ryb:Fund}.

In this paper, we describe the homotopy type 
of the complement $E_\sL$, when $\sL$ is defined over the real numbers,
in terms of the boundary 3-manifold $M_\sL$ of a regular neighborhood 
of $\sL$ in $\C^2$.   We do this by describing $M_\sL$
and $E_\sL$ in terms of a graph of manifolds over the
incidence graph $\Gamma_\sL$.

Let $\Gamma = \{\sV,\sE\}$ be a directed graph, with vertices $\sV$
and edges $\sE$.  Assume that for each oriented 
edge $e \in \sE$ its opposite $\overline e$ is also contained in $\sE$. 
Denote by $i(e)$ the initial point of $e$ and $t(e)$ the terminal point of
$e$.  Thus, $i(e) = t(\overline{e})$ and $t(e) = i(\overline{e})$.
A {\it graph manifold} $\{M_v, M_e, \phi_e\}$ over a directed graph $\Gamma$
is a collection of connected 
vertex manifolds $M_v$, for $v \in \sV$, connected edge manifolds $M_e$, 
with $M_e = M_{\overline e}$, for $e \in \sE$,
and inclusions $\phi_e: M_e \rightarrow M_{t(e)}$ which are isomorphisms
onto a boundary component of $M_{t(e)}$, and which induce endomorphisms on 
fundmantal groups (cf. \cite{Wal:Klasse}).  
A {\it morphism} of graph manifolds over a given
graph $\Gamma$ is a collection
of continuous maps between corresponding vertex and 
edge manifolds which commute with
the maps $\phi_e$.  The manifold associated to a graph manifold is
the space obtained by gluing the $M_v$ together along their boundary
components according to the maps $\phi_e$.  

\begin{theorem}  
If $\sL$ is a complex line arrangement, and $M_\sL$ is the 
boundary 3-manifold of $\sL$, then $M_\sL$ has the structure
of a graph manifold over $\Gamma_\sL$ whose vertex manifolds
are Seifert fibered 3-manifolds.  Furthermore,
if $\sL_1$ and $\sL_2$
are two line arrangements and 
$$
\alpha : \Gamma_{\sL_1} \rightarrow \Gamma_{\sL_2},
$$
is an isomorphism of incidence graphs, then there
is a compatible isomorphism
$$
\beta : M_{\sL_1} \rightarrow M_{\sL_2}
$$
of graph manifolds.
\end{theorem}

Any decomposition of
a 3-manifold by incompressible surfaces gives rise to a presentation
of the manifold as a graph manifold.
Let $H_d$ be a thickened $d$-component oriented Hopf link in $S^3$, 
where each pair of component links having linking number $1$.  Then
$S^3 \setminus H_d$ is a 3-manifold with boundary 2-tori and each
boundary torus has a natural framing $(\mu,\lambda)$.
Theorem 1 is a consequence of the following more explicit description
of $M_\sL$.  The omitted case is treated separately
in the beginning of Section 2.

\begin{theorem} Let $\sL$ be any line arrangement in $\C^2$ so
that each line $L \in \sA$ contains at least one point in $\sP$.  
Then the boundary 3-manifold $M_\sL$ has a torus decomposition
into pieces 
$$
M_v  = S^3 \setminus H_d,
$$
where $v$ ranges over vertices in $\Gamma_\sL$.
If $v$ is a point-vertex, there is a one-to-one correspondence
between the $d$ boundary components of $M_v$ and the 
edges emanating from $v$ in $\Gamma_\sL$.
If $v$ is a line-vertex, there is a one-to-one correspondence 
between $d-1$ of the $d$ boundary components of $M_v$ 
and the edges emanating from $v$ in $\Gamma_\sL$. 
Given an edge $e = e(p,L)$, the corresponding boundary components
$T_{L,p} \subset M_{v_p}$ and $T_{p,L} \subset M_{v_L}$ are identified
by $\phi_e$ and $\phi_{\overline{e}}$ according to the relation:
\begin{eqnarray*}
\mu_{p,L} &=& \lambda_{L,p}\\
\lambda_{p,L} &=& \mu_{L,p} + \lambda_{L,p},
\end{eqnarray*}
where $(\mu_{L,p},\lambda_{L,p})$ is the induced framing of $T_{L,p}$
and $(\mu_{p,L},\lambda_{p,L})$ is the induced framing of $T_{p,L}$.
\end{theorem}

Theorem 2, proved in Section 2, implies that 
the fundamental group of $M_\sL$ is torsion free and residually finite
\cite{Serre:Trees}, \cite{Hemp:Res} (see Corollary 2.4).

Given two topological space $X \subset Y$, let $Y/X$ be the space
$Y$ with $X$ collapsed to a point.
For real line arrangements, the homotopy type of the
complement can be described as follows.

\begin{theorem} Let $\sL$ be a real line arrangement.  Then there 
is a continuous map
$$
f : \Gamma_\sL \rightarrow M_\sL
$$
such that $f(v) \in M_v$ for all vertices $v$ of $\Gamma_\sL$,
and the complement $E_\sL = \C^2 \setminus \sL$ is homotopy
equivalent to $M_\sL/f(\Gamma_\sL)$.
\end{theorem}

The map $f$ of Theorem 3 depends on 
the {\it ordered incidence graph} $\widetilde{\Gamma_\sL}$ 
of $\sL$, the incidence graph $\Gamma_\sL$ 
together with an ordering of the edges emanating from each vertex. 
(An ordered graph is similar to a {\it fat graph}, which has a cyclic
ordering on the edges at each vertex \cite{Pen:Per}.)
A morphism of ordered incidence graphs is a morphism of incidence
graphs which preserves orderings of edges.
Let $\sL$ be a line arrangement defined over the real numbers.
We construct a model for the
homotopy type of $E_\sL$ which only depends on the ordered graph 
$\widetilde {\Gamma_\sL}$ associated to $\sL$ and is simple to describe  
(cf. \cite{Falk:Hom}, \cite{Lib:Hom}, \cite{O-T:Arr},
\cite{Ran:Fund} and \cite{Sal:Top}.)

\begin{theorem}  Let $\sL$ be a real line arrangement.  
Choose basepoints
$b_v \in M_v$ for each vertex $v$ of $\Gamma_\sL$.  
For each edge $e$ in $\Gamma_\sL$, let $\tau_e$ be a path
on $M_{i(e)} \cup M_{t(e)}$ from $b_{i(e)}$ to $b_{t(e)}$
which intersects $M_e$ in one point.  Identify
$\Gamma_\sL$ with its associated singular $1$-complex.  
Then there is a continuous mapping
$$
f: \Gamma_\sL \rightarrow M_\sL
$$
satisfying the following conditions:
\begin{description}
\item{(i)} $f(v)= b_v$ for all vertices $v \in \sV$;
\item{(ii)}  
if $p_1,\dots,p_s$ are the ordered points on $L \cap \sP$,  $e = e(L,p)$,
and $p = p_j$, 
then
$$
f(e) = \left \{\begin{array}{ll}
	\tau_{e}&\qquad\mbox{if $j=1,2$}\\
	g_{L,p_2}\dots g_{L,p_{j-1}}\tau_{e} &\qquad\mbox{if $j>2$},
	\end{array}
	\right .
$$
where $L_1,\dots,L_r$ are the ordered lines through $p$, 
and
$$
g_{L_j,p}= \mu_{L_1,p} \dots \mu_{L_j,p} \mu_{L_{j-1},p}^{-1}
\dots \mu_{L_1,p}^{-1};
$$
\item{(iii)} $f(e) = f(\overline{e})^{-1}$ for all edges $e$; and
\item{(iv)} $E_\sL$ is homotopy equivalent to the space 
$M_\sL/f(\Gamma_\sL)$.
\end{description}
\end{theorem}

Section 3 contains proofs of Theorem 3 and Theorem 4.  In Section 4
we show why the construction fails for general complex line arrangements
and algebraic plane curves.  The author thanks P. Orlik and M. Falk 
for helpful discussions during the writing of this paper.

\section{Boundary manifolds}

Let $\sL$ be a line arrangement in the complex plane $\C^2$, let 
$\sA$ be the set of lines in $\sL$ and let $\sP$ be the points of
intersection on $\sL$.
In this section we describe the boundary manifold of $\sL$
in terms of the incidence graph $\Gamma_\sL$.

\heading{Case of disconnected incidence graph.}
Consider the case when $\sL$ consists of $k$ non-intersecting 
lines.  This is the only case where $\Gamma_\sL$ is not a connected graph
and is a finite union of vertices, one for each line
in $\sA$.   The boundary 3-manifold $M_\sL$ is then
a disjoint union of $k$ solid tori.  
The complement $E_\sL$ is equal to the product 
$$
\C \setminus \{k\ \mbox{points}\} \times \C.
$$
Thus, $M_\sL/f(\Gamma_\sL)$ is homotopy equivalent to $E_\sL$
for any map $f : \Gamma_\sL \rightarrow M_\sL$ which sends
each $v_L$ to a point on the line $L$.

Throughout this paper,
unless otherwise stated we will assume that $\sP$ is non-empty
and therefore $\Gamma_\sL$ is connected.

\heading{Incidence graph.}
The {\it (point/line) incidence graph} $\Gamma_\sL$ of $\sL$ is a 
bipartite graph with {\it point-vertices} 
$$
v_p, \qquad p \in \sP  
$$
and {\it line-vertices} 
$$
v_L, \qquad L \in \sA. 
$$
The edges of $\Gamma_\sL$ are of the form
$$
e(p,L)\ \mathrm{or}\  e(L,p), \qquad p\in \sP, L \in \sA,\ \mathrm{and}\  
p \in L.
$$

The graph $\Gamma_\sL$ is a directed graph.
The {\it initial point}  of $e = e(p,L)$ 
is defined to be $i(e) = v_p$ and the {\it terminal point} 
is defined to be $t(e) = v_L$.  Similarly, 
if $e = e(L,p)$, then $i(e) = v_L$ and $t(e) = v_p$.  We say that 
$e(L,p)$ and $e(p,L)$ are {\it conjugates} of each other and write
$e(L,p) = \overline{e(p,L)}$.

We will also think of $\Gamma_\sL$ as a singular 1-complex 
whose zero cells map to the vertices and whose one cells map
to the edges so that $e$ and $\overline{e}$ are identified but with 
opposite orientations.  The endpoints of each edge $e$ are attached 
to $i(e)$ and $t(e)$ in the obvious way. 

\heading{Graph manifold.}
Let $\Gamma = \{\sV,\sE\}$ be a directed graph such that for each
edge $e \in \sE$, there is a conjugate edge $\overline{e} \in \sE$
so that $i(e) = t(\overline{e})$ and $t(e) = i(\overline{e})$.
A {\it graph manifold}
$\sM$ over $\Gamma$ is a collection $(\{M_v\},\{M_e\},\{\phi_e\})$
where
\begin{description}
\item{(i)} $M_v$ is a connected manifold with boundary 
for all $v \in \sV$;
\item{(ii)} $M_e$ is a connected compact manifold without boundary
for all $e \in \sE$; 
\item{(iii)} $M_e = M_{\overline{e}}$ for all $e \in \sE$; 
\item{(iv)} 
$$
\phi_e : M_e \rightarrow M_{t(e)}
$$
is an embedding of $M_e$ onto a boundary component of $M_{t(e)}$;
\item{(v)} the map $\phi_e$ induces an endomorphism 
$$
{\phi_e}_* : \pi_1(M_e) \rightarrow \pi_1(M_{t(e)})
$$
on fundamental groups.
\end{description}

The underlying space $M$ associated to a graph manifold $\sM$ is
defined to be the space obtained by gluing together the vertex manifolds
$M_v$ along their boundary components so that $\phi_e(M_e)$ is 
identified with $\phi_{\overline{e}}(M_e)$ by 
$\phi_e(q) = \phi_{\overline{e}}(q)$ for all $q \in M_e$.
Any connect sum of manifolds glued along incompressible boundary components
can be thought of as a graph of manifolds.

\heading{Graphs of groups.}
As with graph complexes, one can talk
about the graph of groups associated to a graph manifold.
A graph of groups $\sG$ over a directed graph $\Gamma$ is a collection
of groups $G_v$, for each vertex $v \in \sV$,
groups $G_e = G_{\overline{e}}$, for each edge $e \in \sE$, 
and group endomorphisms 
$$
\psi_e : G_e \rightarrow G_{t(e)}.
$$
The collection of fundamental groups of the vertex
and edge manifolds of a graph of manifolds is a graph of groups.
The underlying group of a graph of groups is obtained from the vertex
and edges groups by a combination of amalgamated products and 
HNN extensions.  If $\sG$ is the graph of groups associated to a graph
of manifolds $\sM$ then the underlying group of $\sG$ is the fundamental
group of the underlying manifold $M$ of $\sM$.

\heading{Regular neighborhood.} 
Consider $\C^2$ as a metric space with the usual
distance function
$$
d((x_1,y_1),(x_2,y_2)) = \sqrt{(x_1 - x_2)^2 + (y_1 - y_2)^2}.
$$
Let $\epsilon > 0$ be such that $d(p,L) > 2 \epsilon$
for all pairs $(p,L) \in \sP \times \sA$ such that $p \notin L$.  

For each $p \in \sP$, let $N_p$ be the ball of radius $\epsilon$ around $p$ in
$\C^2$ and let $S_p$ be the boundary of $N_p$.  Then $N_p$ is a Milnor
ball around $p$ and the pair $(N_p, N_p \cap \sL)$ is a cone over
$(S_p, S_p \cap \sL)$ (see \cite{Milnor:Sing}).  Note that $L \cap N_p \neq \emptyset$ if and
only if $p \in L$.  Let 
$$
\delta_p= \min\{d(L_1\cap S_p, L_2\cap S_p)\} 
$$
where the minimum is taken over $L_1,L_2 \in \sA$ and $p \in L_1 \cap L_2$.
Take $\delta > 0$ such that $\delta_p > 2\delta$ for all $p \in \sP$.
For each line $L \in \sA$, let 
$$
N_L  = \{q \in \C^2\ : \ d(q,L) \leq \delta\}.
$$
By choosing $\delta$ smaller if necessary, we can assume that for each
edge $e(p,L)$ in $\Gamma_\sL$, $N_L$ meets $S_p$ transversally.
Let 
$$
N_\sL = \bigcup_{p \in \sP} N_p \cup \bigcup_{L \in \sA} N_L.
$$
We call $N_\sL$ a {\it regular neighborhood} of $\sL$ in $\C^2$.

\heading{Boundary manifold.}
Consider the projective compactification $\P^2$ of $\C^2$ and let
$L_\infty$ be the line at infinity, so that $\P^2 = \C^2 \cup L_\infty$.
Let $Q$ be a point on $L_\infty$ not on the projective closure 
$\overline{\sL}$ of $\sL$.  Let $\overline L$ be the projective closure
of $L$ for all $L \in \sL$.  

For each $L \in \sA$, let $B_L$ be the closure
of $N_L$ in $\P^2$.  Then $B_L$ is a closed tubular neighborhood of 
$\overline L$ in $\P^2$. 
Let $B_\infty$ be a closed tubular neighorhood of $L_\infty$ not containing
any element of $\sP$.  
Let $S_L, S_\infty$ be the boundaries of
$B_L$ and $B_\infty$, respectively.

The boundary $M_\sL$ of $N_\sL$ has a deformation retraction onto
$M_\sL \setminus B_\infty$.
Hereafter, for simplicity, we will replace $M_\sL$
by $M_\sL \setminus B_\infty$.

For each $p \in \sP$, note that
$$
\bigcup_{L \in \sA} N_L
$$
intersects $S_p$ transversally in disjoint 2-tori.
Set 
$$
M_p = S_p \setminus \bigcup_{L \in \sA} N_L.
$$
For each $L \in \sA$, let
$$
M_L = S_L \setminus (B_\infty \cup \bigcup_{p \in \sP} N_p).
$$
Then $M_\sL$ is a connect sum of the $M_p$ and $M_L$. 

By our construction we have the following.

\begin{lemma} 
The submanifolds $M_p$ and $M_L$ in $M_\sL$ intersect in a boundary 
component if and only if $e(p,L)$ is an edge in $\Gamma_\sL$, and 
this intersection is a $2$-dimensional torus.
\end{lemma}

Recall that a Hopf fibration $h : S^3 \rightarrow S^2$ is an oriented 
circle bundle such that the fibers have intersection number 1.  Given any
$d$ points $\sQ$ in $S^2$, the preimage $h^{-1}(\sQ)$ is called a
$d$-component Hopf link (see Figure \ref{fig:Hopf}).  The 
pair $(S^3, h^{-1}(\sQ))$ does not
depend on which $d$ points are chosen.  Let $H_d$ be a thickening
of $h^{-1}(\sQ)$.  

\makefig{3-component Hopf link.}{fig:Hopf}
{\psfig{figure=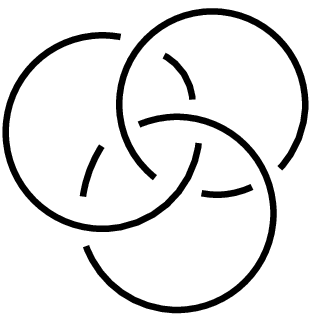,height=1in}}

\begin{lemma} 
For each $p \in \sP$ there is a natural identification of 
the pair $(S_p,M_p)$ with $(S^3,S^3 \setminus H_d)$, 
where $d$ is the degree of $v_p$.
For each $L \in \sA$, there is a natural identification of
the pair $(S_L,M_L)$ with the pair $(S^3,S^3 \setminus H_{d+1})$, where $d$
is the degree of $v_L$.
\end{lemma}

\heading{Proof.} 
We begin with $M_p$.  Identify $S_p$ with the 3-sphere $S^3$ and
$L_\infty$ with $S^2$.  Then the pencil of lines through $p$ determines
a projection $\rho: S_p \rightarrow L_\infty$ giving a commutative diagram
$$
\cd
{
&S_p &\mapright{=} &S^3\cr
&\mapdown{\rho} &&\mapdown{h}\cr
&L_\infty &\mapright{=} &S^2
} 
$$
where $h$ is the standard Hopf fibration.  

\makefig{Neighborhood of a point $p \in \sP$.}{fig:pnbd}
{\psfig{figure=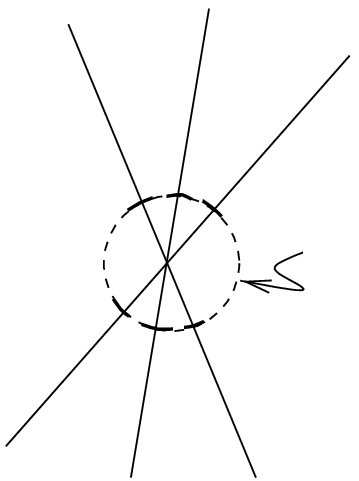,height=1.5in}}

Let $L_1,\dots,L_d \subset \sL$ be the lines passing through $p$.
Under the identification of $S_p$ with $S^3$, $M_p$ is identified with
the complement in $S^3$ of 
$$
h^{-1}(\bigcup_{i=1}^d D_i),
$$
where each $D_i$ is a small disk around the point in $S^2$ corresponding 
to $\overline{L_i} \cap L_\infty$.
We thus have a natural identification of the pair 
$(S_p,M_p)$ with $(S^3,S^3 \setminus H_d)$.


Describing $M_L$ is for the most part the dual picture to the above.
Let $S_L$ be the boundary of a tubular neighborhood of $S_{\overline L}$.
Let $B_Q$ be a ball neighborhood of $Q$ in $\P^2$ and let $S_Q$ be 
its boundary 3-sphere.  Then $S_L$
and $S_Q$ are canonically identified as 3-spheres fibering over a
2-sphere.
The pencil of lines through $Q$ defines a commutative diagram
$$
\cd
{
&S_Q &\mapright{=} &S^3 \cr
&\mapdown{\rho} && \mapdown{h}\cr
&\overline L &\mapright{=} &S^2
}
$$
where again $h$ is the standard Hopf fibration.
Under the identification of $S_Q$ with $S_L$, $B_\infty$
equals $h^{-1}(D_\infty)$ where $D_\infty$ is a small disk neighborhood
of $\overline L \cap L_\infty$.

\makefig{Neighborhood of a line $L \in \sA$.}{fig:lnbd}
{\psfig{figure=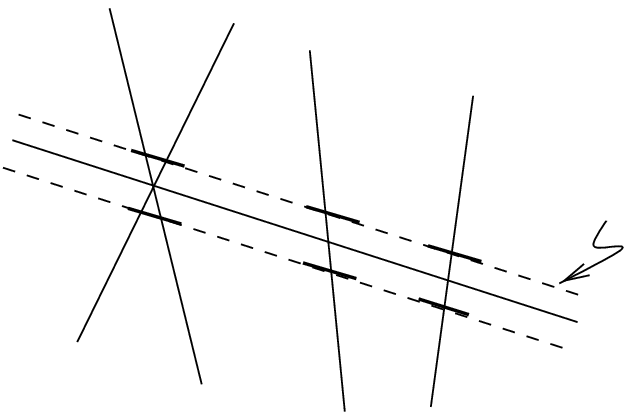,height=1.5in}}

Let $p_1,\dots,p_d \in \sP$ be the points in $\sP$ which lie on $L$.
Then the identification of $S_L$ with $S_Q$ gives 
an identification of $M_L$ with
$$
S_Q \setminus h^{-1}(D_\infty \cup \bigcup_{i=1}^d D_i),
$$
where each $D_i$ is a small disk neighborhood of $p_i$ in $\overline L$.
This identifies $(S_L,M_L)$ with $(S^3, S^3 \setminus H_{d+1})$. \qed

Lemma 2.1 and Lemma 2.2 imply the following.

\begin{proposition} For any line arrangement $\sL$, $M_\sL$ is a Haken
3-manifold with a torus decomposition into Seifert fibered manifolds.
This torus decomposition gives $M_\sL$ the structure of a graph manifold
over $\Gamma_\sL$.
\end{proposition}

\heading{Proof.}  We only need to check when the torus boundary components 
of $S^3 \setminus H_d$ are incompressible.  This is true as long as 
$d > 1$, which holds unless $\sL$ is
a union of parallel lines.\qed

\begin{corollary} The fundamental group $\pi_1(M_\sL)$ is torsion
free and residually finite.
\end{corollary}

\heading{Proof.}
The fundamental group of $\pi_1(M_\sL)$ is a graph of groups whose
vertex groups are $\pi_1(S^3 \setminus H_d) = \Z \times F_d$, where
$F_d$ is a free group on $d-1$ generators, and are torsion free.
It follows that $\pi_1(M_\sL)$ is also torsion free 
(see, for example, Section 1.3, 
Corollary 2 and Section 1.4, Proposition 5 of \cite{Serre:Trees}.) 

Since $M_\sL$ is a 3-manifold with 
a torus decomposition along incompressible tori
into Seifert fibered pieces, by
\cite{Hemp:Res}, Theorem 1.1, the fundamental
is residually finite.  
\qed  

\heading{The gluing maps.} 

Let $T$ be a 2-dimensional torus.  A {\it framing} $(\mu,\lambda)$
on $T$ is a choice of generators for its fundamental group $\pi_1(T)$.

As stated in Proposition 2.3, the edge manifolds of $M_\sL$ considered
as a graph manifold are all 2-dimensional tori, so they are isomorphic 
to $T$.  To prove Theorem 2, we need to describe the gluing maps 
$$
\phi_e : T \rightarrow M_p
$$
and
$$
\phi_{\overline{e}} : T \rightarrow M_L
$$
for each edge $e = e(p,L)$ on the graph $\Gamma_\sL$. 

By Lemma 2.2, the vertex manifolds $M_p$ and $M_L$ are naturally 
identified with 
$E_d = S^3 \setminus H_d$ for some integer $d \ge 2$.
Furthermore, the complex structure on $M_p$ and $M_L$ 
determines an orientation on the core loops of $H_d$.
The oriented pair $(S^3,H_d)$ determines
a framing on the boundary components of $E_d$ and hence on
the boundary components of $M_p$ and $M_L$.  

Let $T$ be a boundary component of $E_d$ and let $\ell$ be
its core curve in $S^3$.  The framing $(\mu,\lambda)$ on
$T$ is given as follows.  Choose a basepoint $x$ on $T$. 
The loop $\mu$ is a positively oriented meridian loop around the 
core curve based at $x$.  If one considers
$T$ as an $S^1$-bundle over $\ell$, then $\mu$ is a loop going
once in the positive direction around the fiber of the bundle.  
The loop $\lambda$ is
a loop based at $x$ ``parallel" to the core curve.  That is, it 
is a positively
oriented loop whose linking number with the core curve is zero.
These definitions uniquely determine $\mu$ and $\lambda$ up
to homotopy.   

Let $p \in \sP$, $L \in \sA$ be such that $p \in L$.
We will write $T_{L,p}$ for the boundary component on 
$M_p$ corresponding to $L$, and
$T_{p,L}$ for the boundary component on $M_L$ corresponding to $p$.
Let $(\mu_{L,p}, \lambda_{L,p})$ be the 
framing on $T_{L,p}$ and $(\mu_{p,L},\lambda_{p,L})$ be the framing
on $T_{p,L}$.

To describe $\phi_e$ up to homotopy, it suffices to describe what
it does to the framings.  

\begin{lemma}  Let $e = e(L,p)$ be an edge of $\Gamma_\sL$.  If 
we identify $T$
with $T_{L,p}$ by the map
\begin{eqnarray*}
\phi_e : T &\rightarrow& T_{L,p} \subset M_p\\
\phi_e(\mu) &=&  \mu_{L,p}\\
\phi_e(\lambda) &=& \lambda_{L,p},
\end{eqnarray*}
then
\begin{eqnarray*}
\phi_{\overline{e}} : T &\rightarrow& T_{p,L} \subset M_L\\
\phi_{\overline{e}}(\mu) &=& \mu_{p,L} + \lambda_{p,L}\\
\phi_{\overline{e}}(\lambda) &=& \mu_{p,L}
\end{eqnarray*}
gives the gluing map for the graph manifold $M_\sL$.
\end{lemma}

\makefig{Intersection of $S_L$ and $S_p$.}{fig:gluing1}
{\psfig{figure=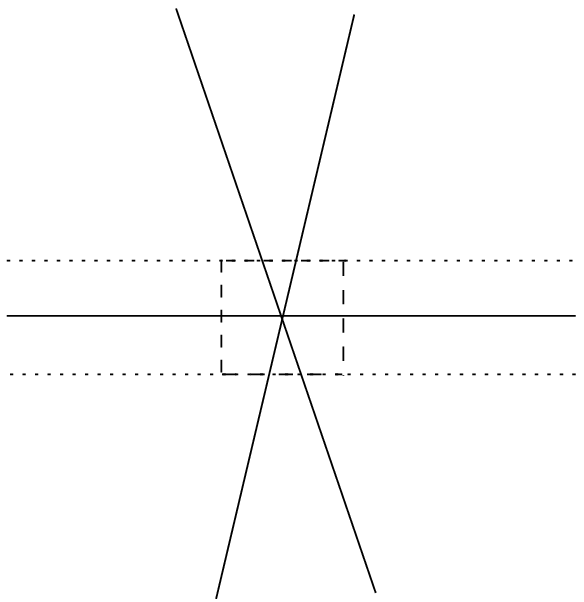,height=1.5in}}

\heading{Proof.}  We will study how $S_L$ and $S_p$
meet by changing coordinates so that $L$ is given by the
equation $y=0$, the point $p$ is the origin, $p = (0,0)$. 
Then $S_L$ and $S_p$ are given by
$$
S_L = \overline{\{(x,y) \in \C^2\ : \ |y| = 1\}},
$$
and 
\begin{eqnarray*}
S_p &=& 
\{(x,y) \in \C^2\ : \ |x| \leq 1, |y| = 1\}
\cup 
\{(x,y) \in \C^2\ : \ |x| = 1, |y| \leq 1\}\\
&=& \T_1 \cup \T_2
\end{eqnarray*}
as shown in Figure \ref{fig:gluing1}.  
Then $S_L \cap S_p = T$, where $T$ is the 2-torus
$$
T = \{(x,y) \in \C^2\ : \ |x|=1, |y|=1\}
$$
and is the common boundary of $\T_1$ and $\T_2$.

As before let $T_{L,p}$ be the boundary component of $S_p$
corresponding to $L$ and $T_{p,L}$ be the boundary component
of $S_L$ corresponding to $p$.  Since the components of the
Hopf link are all unknots, the framing on the link boundary
components is given as follows.  Each link component has
a solid torus neighborhood $N$ in $S^3$. The complement of $N$
is also a solid torus $N^c$.  The framing on the boundary of $N$
is given by $(\mu,\lambda)$ where $\mu$ contracts in $N$ and
$\lambda$ contracts in $N^c$. 

Thus, $T_{L,p}$ has the framing
\begin{eqnarray*}
\mu_{L,p} &=& (1, \exp(2\pi i t))\\
\lambda_{L,p} &=& (\exp(2\pi i t),1)
\end{eqnarray*}
and $T_{p,L}$ has the framing
\begin{eqnarray*}
\mu_{p,L} &=& (\exp(2\pi i t), 1)\\
\lambda_{p,L} &=& (\exp(2 \pi i t), \exp (2 \pi i t)).
\end{eqnarray*}

This implies $\mu_{p,L} = \lambda_{L,p}$ and
$\lambda_{p,L} = \mu_{L,p} + \lambda_{L,p}$,
as illustrated in Figure \ref{fig:gluing2}. 
\qed

\makefig{Gluing map.}{fig:gluing2}
{\psfig{figure=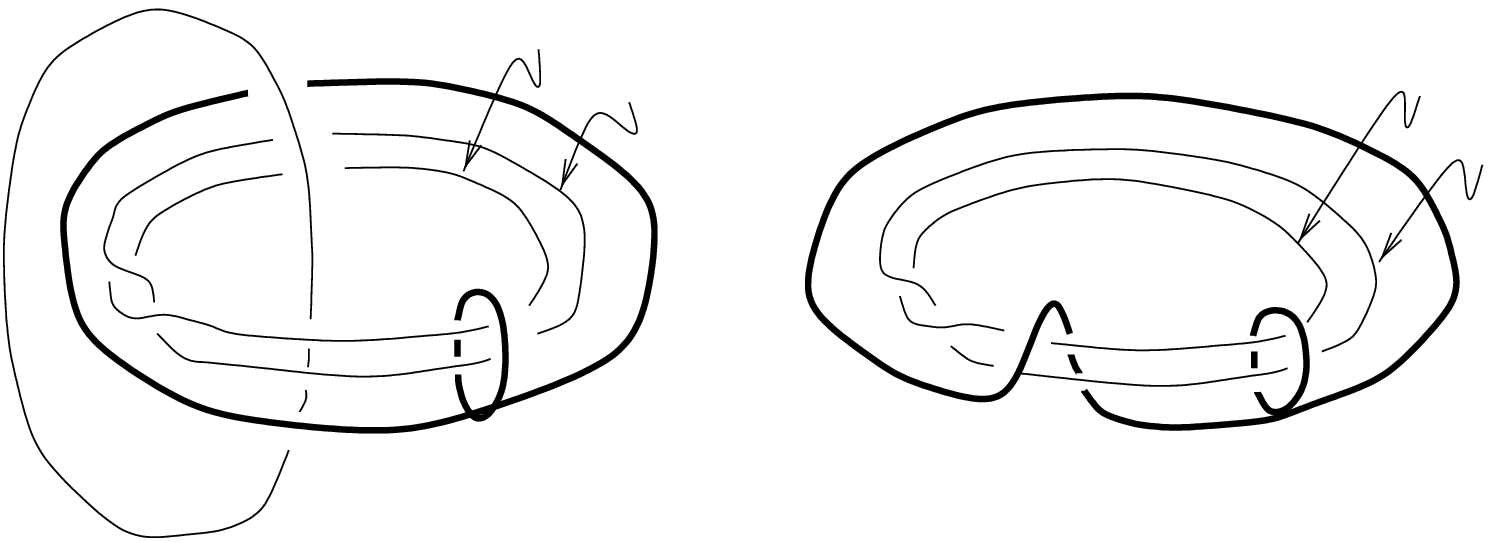,height=1.75in}}

\heading{Proof of Theorem 2.} 

We have shown that $M_\sL$ is a graph manifold whose vertex manifolds
are $S^3 \setminus H_d$, where $d$ is the degree of the vertex if
the vertex is a point-vertex and $d-1$ is the degree of the vertex 
if the vertex is a line-vertex.  The identifications of
the edge manifolds with the boundary components of the vertex manifolds
can be described in terms of the natural framings on the boundary 
components of $S^3 \setminus H_d$ by the edge maps described in Lemma 2.5.
This completes the proof of Theorem 2.\qed

\section{Real line arrangements}

In this section, we concentrate on real line arrangements, and will assume
that $\sL$ is defined by equations of the form
$$
y=a_i x + b_i, \qquad i=1,\dots,k,
$$
where $a_i,b_i \in \R$.  We will show how to reconstruct the homotopy
type of the complement of $\sL$ in terms of the ordered graph of $\sL$.

\heading{Ordered graphs.}

\makefig{Ordered graph associated to the Ceva arrangement.}{fig:graph}
{\psfig{figure=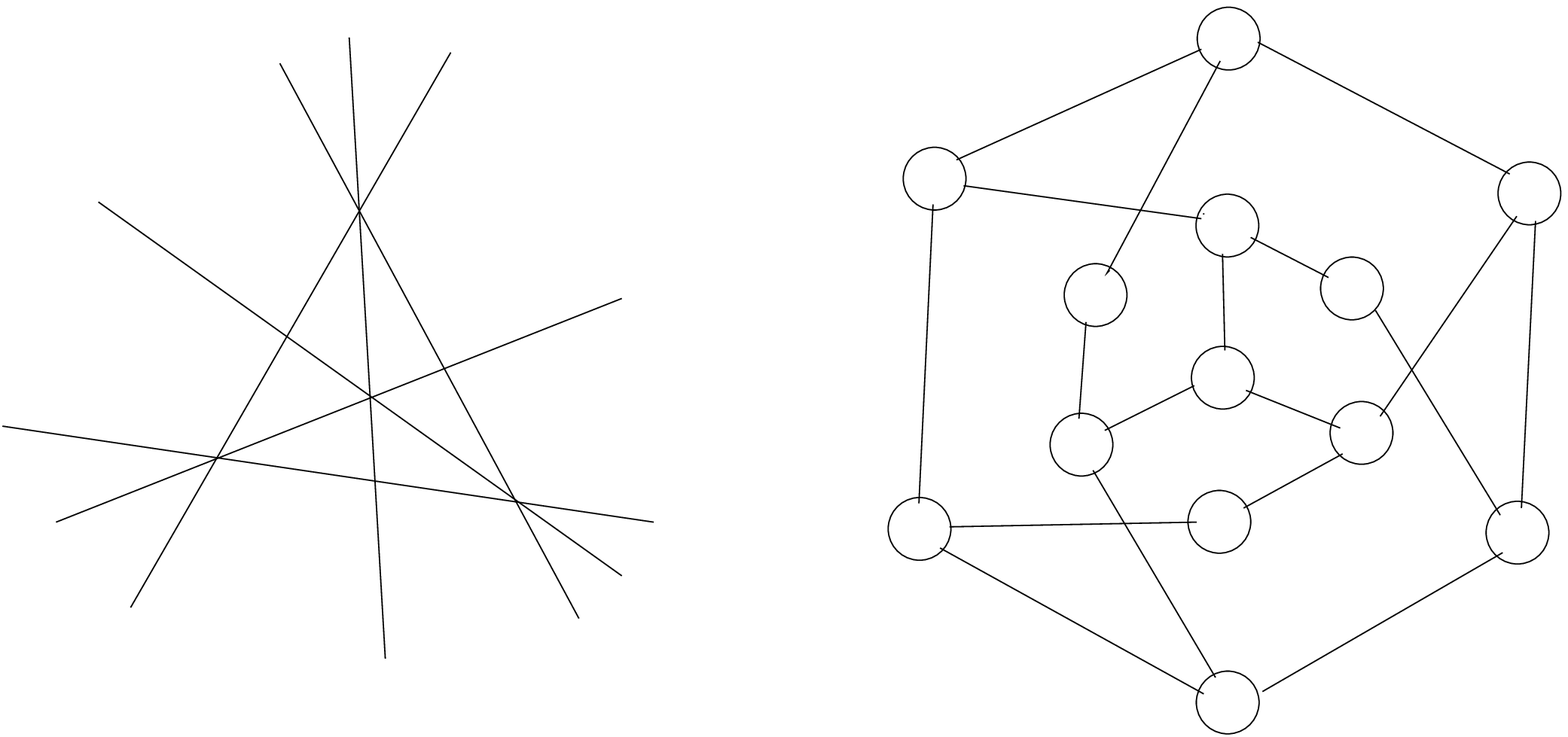,height=2in}}

The ordered graph associated to $\sL$ is the incidence graph $\Gamma_\sL$
together with some extra structure.
For real line arrangements, we will order the 
edges emanating from each vertex of the incidence graph $\Gamma_\sL$
as follows.  If $v_p$ is the point-vertex associated to the point 
$p \in \sP$, then we order the edges $e(p,L_1),\dots,e(p,L_r)$ emanating from
$p$ by the slopes of $L_1,\dots,L_r$ in decreasing order.
If $v_L$ is the line-vertex associated to the line $L \in \sA$, 
we order the edges $e(L,p_1),\dots,e(L,p_s)$ emanating from $v_L$
by the $x$-coordinates of $p_1,\dots,p_s$ in decreasing order.
Since we assume that
none of the lines are parallel to the $y$-axis, this is well-defined.
The incidence graph $\Gamma_\sL$ of $\sL$ endowed with these orderings
on the edges emanating from vertices is called the {\it ordered graph} 
associated to $\sL$.  For example, Figure \ref{fig:graph} gives the 
ordered graph of the Ceva arrangement.   

The global ordering of the points in $\sP$ by their
$x$-coordinate and the lines in $\sA$ by their slope determines
the ordered graph.  The orderings near $v_{L_1}$ and $v_{p_5}$
are given are given in Figure \ref{fig:graph}.

\heading{The skeleton of a line arrangement.}

Let $M_\sL$ be the boundary manifold of $\sL$ and let
$$
\alpha: M_\sL \rightarrow \sL
$$
be the natural projection map, well-defined up to homotopy.
Let $\Sigma_\sL$ be the
real part of $\sL$ with all infinite ends removed.  We will call 
$\Sigma_\sL$ the skeleton of $\sL$.  
Figure \ref{fig:skeleton} gives the skeleton of the Ceva arrangement.

\makefig{Skeleton associated to the Ceva arrangement.}{fig:skeleton}
{\psfig{figure=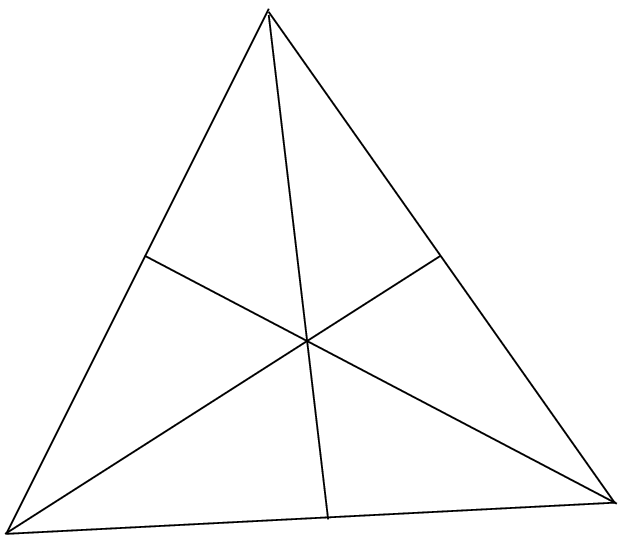,height=1in}}

We start by describing the homotopy type of $E_\sL$ in terms of 
$M_\sL$ and $\Sigma_\sL$.

\begin{lemma} There is an embedding 
$$
g : \Sigma_\sL \rightarrow M_\sL
$$
so that 
\begin{description}
\item{(i)} $\alpha\circ g$ is the identity map; 
\item{(ii)} the image of $g$ contracts in $E_\sL$; and
\item{(iii)} $E_\sL = M_\sL/g(\Sigma_\sL)$.
\end{description}
\end{lemma}

\heading{Proof.}
The embedding $g$ is given as follows.  
If $(x,y) \in N_p$, for some $p=(x_p,y_p) \in \sP$, then 
$$
g(x,y) =  (x,y + i \max\{{\sqrt{\epsilon^2-(y-y_p)^2},\delta}\}),
$$
otherwise
$$
g(x,y) = (x,y+i\delta).
$$
The image $g(\Sigma_\sL)$ contracts in $E_\sL$, since in $E_\sL$
it is isotopic to $\Sigma_\sL + (0,i\epsilon) \subset \R^2 + (0,i\epsilon)$.
No point $(x,y+i\epsilon) \in \R^2 + (0,i\epsilon)$ satisfies
an equation of the form
$$
y=ax +b 
$$
where $a$ and $b$ are real numbers.
Thus, $\R^2 + (0,i\epsilon)$ is contained in $E_\sL$.  Since
$\R^2 + (0,i\epsilon)$ is contractable, we see that $g(\Sigma_\sL)$
is also contractable.

Now let us consider the projection of $\rho : \C^2 \rightarrow \C$
onto the first coordinate. Let $\rho_\sL$ be the restriction of $\rho$
to $E_\sL$.  Then $\rho_\sL$ is a topological fibration over the 
complement of $\sQ = \rho(\sP)$.  

\makefig{Contraction of the $x$-coordinate plane to $V_Q$.}{fig:contraction}
{\psfig{figure=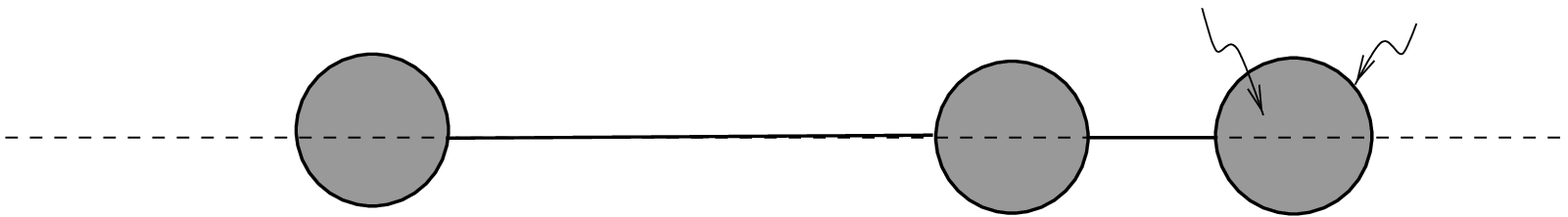,width=4in}}

For each $q \in \sQ$ let 
$$
D_q = \{x \in \C\ : \ |x - q|\leq 1\}
$$
and let $S_q$ be the boundary of $D_q$.
Assume (by expanding coordinates if necessary) that 
the $S_q$ do not intersect one another.  Let $W_Q$ be the
union of the $S_q$ and real line segments joining the $S_q$,
and let $V_Q$ be the union of $W_Q$ and the $D_Q$
as in Figure \ref{fig:contraction}.

The deformation retraction of $\C$ onto 
$V_Q$ extends to a deformation retraction of $E_\sL$ onto
$\rho_\sL^{-1}(V_Q)$.

The fibers of $\rho_\sL$ over $V_{\sQ}$ split up into 3 types: those
over the line segments in $W_{\sQ}$; those over $S_q$; and those 
over the $q$.  These retract to the outlined and shaded regions
shown in Figure \ref{fig:fibers} which we will write as
$F_x$, for $x \in I$,
$F_s$, for $s \in S$ and $F_q$, respectively. 

\makefig{Fibers over $V_{\sQ}$.}{fig:fibers}
{\bigskip\psfig{figure=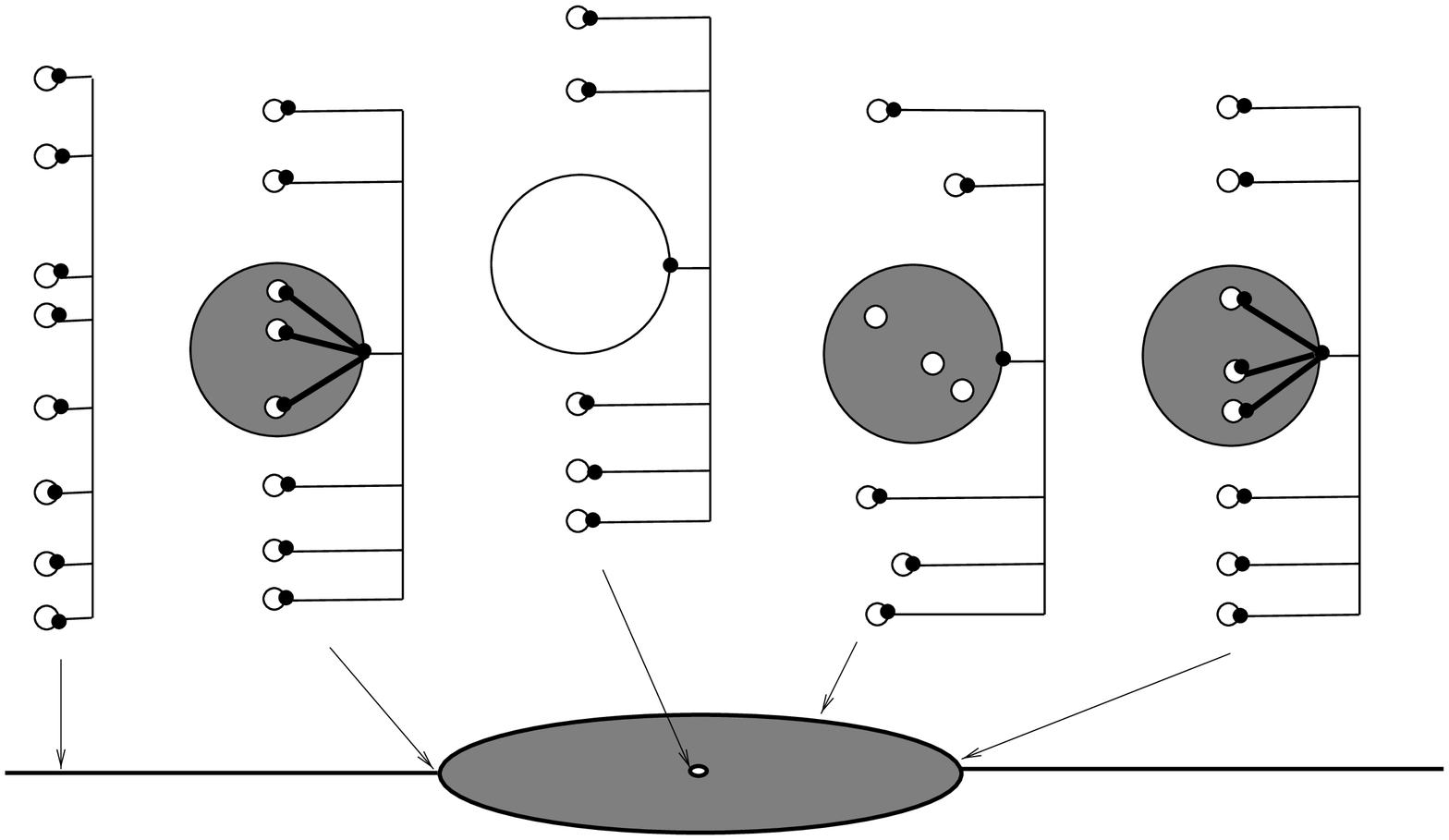,height=2.75in}}

We reconstruct the homotopy type of $E_\sL$ from this picture 
as follows.   First notice that the space $\rho^{-1}|_{M_\sL}(D_q)$ has
a deformation retraction to 
$$
F_{D_q} = D_q \times F_q \cup \bigcup_{s \in S_q}F_s. 
$$
Under this retraction $g(\Sigma_\sL)$ maps to the right most points of
the circles in $F_s$, $F_x$ and 
$F_q$ union the line segments depicted on
$F_{s_1}$ and $F_{s_2}$.  
The union of the right-most points of the circles on
$F_s$ as $s$ ranges in $S_q$ bounds a $2$-disk 
$$
\{p + (0,iy): |y| \leq \epsilon\} = D_q \times \{p + (0,i\epsilon)\}
$$
in the deformation of $\rho^{-1}|_\sL(D_q)$.
Thus, if we include all right-most points of circles in the fibers over
$F_x,F_s,$ and $F_q$, union the line segments on
$F_{s_1}$ and $F_{s_2}$, we obtain a set $G$ which retracts onto 
$g(\Sigma_\sL)$ in $E_\sL$.  
\qed

\heading{Proof of Theorem 4.}

The existence of a map $f : \Gamma_\sL \rightarrow M_\sL$
with property $(iii)$ of Theorem 4 follows from Lemma 3.1, since
$\Gamma_\sL$ can be continuously deformed to $\Sigma_\sL$.  One can
see this by noting that the figures in Figure \ref{fig:compare}
are homotopy equivalent.

\makefig{Incidence graph and skeleton.}{fig:compare}
{\psfig{figure=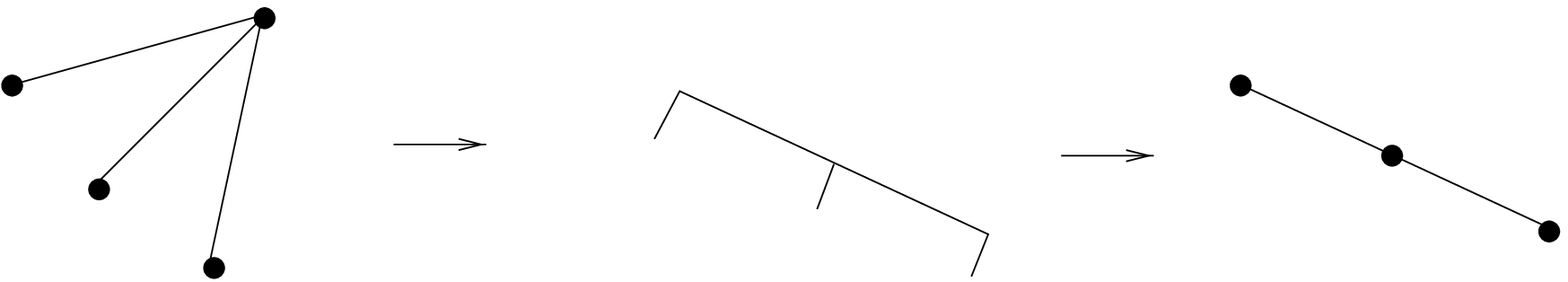,height=1in}}

We will explicitly describe the homotopy type of a map $f : \Gamma_\sL
\rightarrow M_\sL$ satisfying the conditions of Theorem 4
in terms of the ordered graph.

We begin by defining a map $\sigma : \Gamma_\sL \rightarrow \Sigma_\sL$. 
For the point-vertices, we simply send $v_p$ to the point $p$ on 
$\Sigma_\sL$.  For the line vertices, there is no canonical choice.

Let $L \in \sA$.  If $\degree(v_L) = 1$, then let $\sigma(v_L) = p$
where $p = \sP \cap L$.  Otherwise, let $p_1,p_2,\dots,p_r$ be the
points in $\sP \cap L$.  Choose a point on  
the line segment on $\Sigma_\sL$ strictly between $p_1$ and $p_2$
and map $v_L$ to this point. 

For the edges of $\Gamma_\sL$, once we have determined where the
vertices go, the edges map to the unique straight line segment
on $\Sigma_\sL$ connecting the endpoints.  This defines a continuous
map $\sigma: \Gamma_\sL \rightarrow \Sigma_\sL$. 
The composition $f = g \circ \sigma$ gives a map satisfying
properties (i),(iii) and (iv) of Theorem 4.

For each vertex $v$ on $\Gamma$, let $z_v = f(v)$.  (We will write
$z_p = z_{v_p}$ and $z_L = z_{v_L}$.)
To describe the homotopy type of the map $f$, it suffices to describe
the homotopy type of $f(e)$ for each edge $e$ in $\Gamma$.
The path $\sigma(e)$ breaks up into segments
each passing over one point $p \in \sP$.   
Let $I_{L,p}$ be a straight line segment on $\Sigma_\sL \cap L$ such 
that $I_{L,p} \cap \sP = p$.  We will describe the homotopy type
of $g(I_{L,p})$.

\makefig{Path $\gamma_L$ lifting to $L_0$.}{fig:gamma}
{\psfig{figure=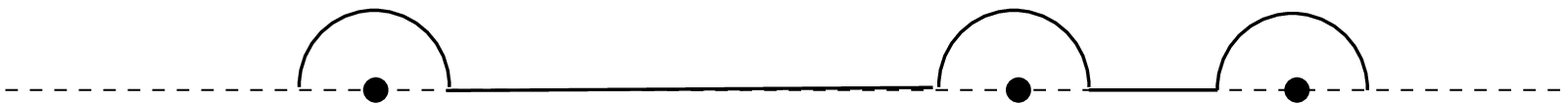,width=4in}}

First fix $L$.
Consider the path $\gamma_L$ on the complex plane
shown in Figure \ref{fig:gamma}.
Here $\sQ = \{q_1,\dots,q_s\} = \rho(\sP \cap L)$.
The projection map $\rho$ is one-to-one when restricted to $L$.  Since
$\rho^{-1}(\gamma_L) \cap L$ is contained in $L_0$, $b_L = \rho^{-1}(\gamma_L)
\cap L$
defines a contractable subset of $L_0$.  
Let $h : \gamma_L \rightarrow M_L$ be defined by 
$$
h(x,y) = \rho|_L^{-1}(x,y) + (0,\delta i),
$$
for each $(x,y) \in \gamma_L$. 
Then $h(\gamma_L)$ defines a contractible subset of $M_L$.
We can assume that the image contains the basepoint 
$z_L$.

\makefig{Generators for $M_p$.}{fig:generators}
{\psfig{figure=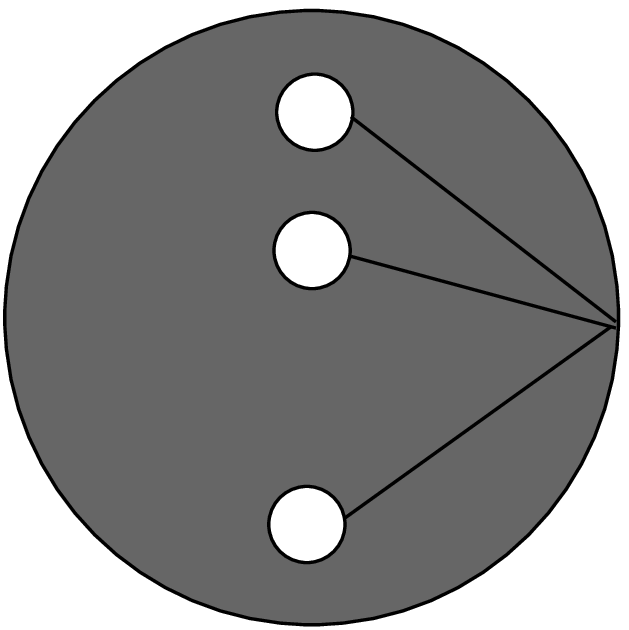,height=0.75in}}

Now fix $p \in L \cap \sP$.
Let $q = \rho(p)$ and let $J_p$ be the
arc segment of $\gamma_L$ near $q$. 
For each fiber $F_s$ over points $s \in J_p$, $h(J_p) \cap F$ 
is the rightmost point in the circle (see Figure \ref{fig:homotopy})
on $F$ corresponding to $L \cap F$.
As noted in the proof of Lemma 3.1, the 
path $g(I_p)$ is isotopic to the path whose intersection with
each fiber $F_s$ over an interior point $s$ of $J_p$ 
is the right-most point of the large circle
and whose intersection with the fibers $F_{s_1}$ and $F_{s_2}$
is a path
from the right-most point of the large circle to the right-most
point of the inner circle associated to $L$ (see Figure {fig:fibers}).
The large disk in the fiber $F_s$
over $s \in J_p$ rotates in the counter-clockwise
direction by 180 degrees as $s$ moves from right to left on $J_p$,
as illustrated in Figure \ref{fig:homotopy}.

\makefig{Homotopy type of $f(e(p,L))$.}{fig:homotopy}
{\psfig{figure=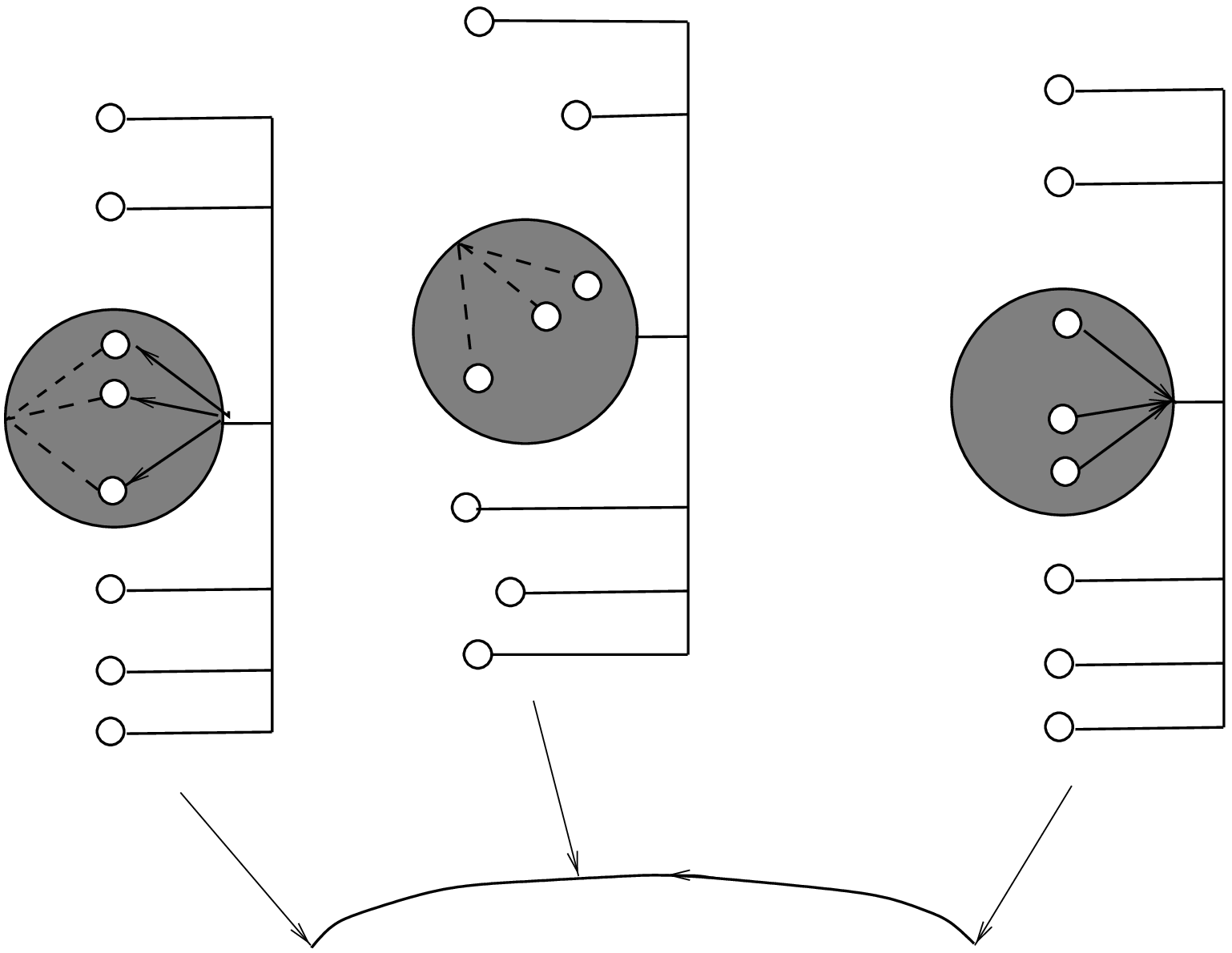,height=2.5in}}

Note that the basepoint $z_p$ for $M_p$ is the right-most point on 
the large circle on $F_{s_1}$ and that $b_L \cap F_s$ 
is the right-most point on the inner circle corresponding
to $L$ for all $s \in J_p$.  
Let $\tau_{L,p}$ be the straight line segment from 
$b_L \cap F_{s_1}$ to $z_p$ (see Figure \ref{fig:homotopy}.)

Let
$\mu_{L_1,p},\dots,\mu_{L_r,p}$ be the generators for $\pi_1(M_p)$ 
as in Figure \ref{fig:generators}.
Then $g(I_p)h(J_p)^{-1}$ is homotopy equivalent to the path $abc$
on $F_{s_1}$ pictured in Figure \ref{fig:paths}.  
In this example $L_1,L_2,L_3$ 
are the lines passing through $p$ in the order of their slopes and
$L = L_2$, then $h(J_p)$ is homotopy equivalent to 
$g(I_p)h(J_p)^{-1} = abc = \tau_{L,p}\mu_{L_3,p}^{-1}\mu_{L_2,p}\mu_{L_3,p}
\tau_{L,p}^{-1}$.  

\makefig{Homotopy type of $g(I_p)h(J_p)^{-1}$}{fig:paths}
{\psfig{figure=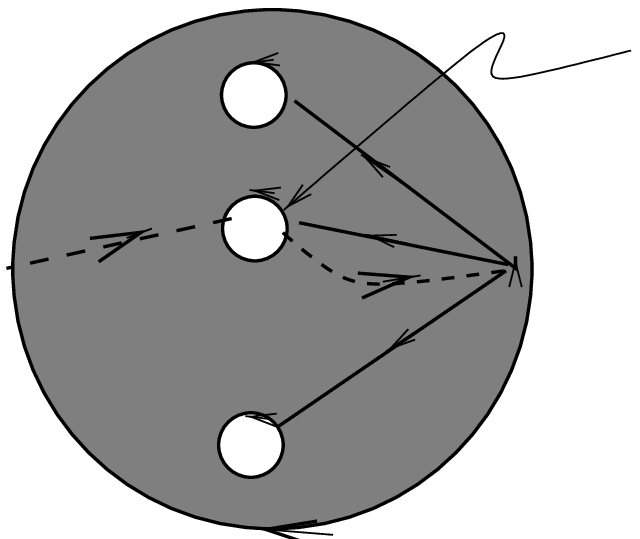,height=1.25in}}

In the general case, we have the following lemma.

\begin{lemma} Let $p \in \sP$ and let $I_p$ be a line segment
on $\Sigma_\sL \cap L$ such that $I_p \cap \sP = \{p\}$
and $p$ lies in the interior of $I_p$.
Let $L_1,\dots,L_r$ be an ordering by slope of the lines in 
$\sA$ passing through $p$.  Let $j$ be such that $L_j = L$.
Then the homotopy class of the difference $g(I_p)h(J_p)^{-1}$ on 
$M_\sL$ is given by the element 
$$
g_{L,p}=\tau_{L,p}\mu_{L_r,p}^{-1}\mu_{L_{r-1},p}^{-1}\dots \mu_{L_j,p}
\mu_{L_{j+1},p} 
\dots \mu_{L_r,p}\tau_{L,p}^{-1}.
$$ 
\end{lemma}

We now describe the homotopy type of the map
$$
f : \Gamma_\sL \rightarrow M_\sL
$$
in terms of the ordered graph $\Gamma_\sL$.

Fix $p \in \sP$ and 
let $e(p,L_1),\dots,e(p,L_r)$ be the ordered edges emanating from $v_p$. 
For each line $L \in \sA$, let $e(L,p_1),\dots,e(L,p_s)$ be the
ordered edges emanating from $v_L$.

For any edge $e(L,p)$ in $\Gamma_\sL$, replace $\tau_{L,p}$ by the
path from $z_L \in M_L$ to $z_p \in M_p$ which goes along
$h(\gamma_L)$ and the original $\tau_{L,p}$.
This does not change the homotopy type of $\tau_{L,p}$.
Note that the choices of $\tau_{L,p}$, $z_p$, $z_L$
a continuous map
$$
\ell : \Gamma_\sL \rightarrow M_\sL
$$
where $\ell(v) = z_v$, for vertices $v$  and
$\ell(e(L,p)) = \tau_{L,p}$ for edges $e(L,p)$.
We will call this map $\ell$ the {\it lifting}
of $\Gamma_\sL$ in $M_\sL$ defined by the collection
$\tau_{L,p}$.
Let $g_{L,p} \in \pi_1(M_p,z_p)$ be defined as in Lemma 3.2,
except that its basepoint is at $z_L$.
We have shown the following.

\begin{lemma}   Let 
$$
f : \Gamma_\sL \rightarrow M_\sL
$$
be the continuous map on $\Gamma_\sL$ considered as a $1$-complex
defined by the following data:  
\begin{description}
\item{(i)}  for the point-vertex $v_p$, $f(v_p) = z_p$;
\item{(ii)} for the line-vertex $v_L$, $f(v_L) = z_L$; and
\item{(iii)} for the edge $e(L,p)$, if
$e(L,p_1),\dots,e(L,p_s)$ are the ordered edges emanating 
from $v_L$, and $p = p_j$, then 
$$
f(e(L,p)) 
= \left \{ \begin{array}{ll}
\tau_{L,p}&\qquad \hbox{if $j = 1,2$;}\\
g_{L,p_2}\dots g_{L,p_{j-1}} \tau_{L,p} &\qquad \hbox{if $j>2$}
\end{array}
\right .
$$
\end{description}
and  $f(e(p,L)) = f(e(L,p))^{-1}$.
Then $E_\sL$ is homotopy equivalent to $M_\sL / f(\Gamma_\sL)$.
\end{lemma}

For each pair $(p_0,L_0) \in \sP \times \sA$ with $p_0 \in L_0$,
the element $g_{L_0,p_0}$ depends on the ordering of 
the edges $e(p_0,L)$ emanating from $v_{p_0}$ and $f(e(L_0,p_0))$ 
depends on the ordering the edges $e(L_0,p)$ emanating from
$v_{L_0}$.  This completes the proof of Theorem 4.\qed

\heading{Fundamental group.}

Choose a vertex $v$ of $\Gamma_\sL$.  We present the fundamental
group $G_\sL = \pi_1(M_\sL,f(v))$ in the notation of \cite{Serre:Trees}
(p. 41).
For each edge $e = e(L,p)$ in $\Gamma_\sL$, let 
$\tau_e = \tau_{L,p}$ and let $\tau_{\overline e} = \tau_{e}^{-1}$.
The corresponding lifting of $\Gamma_\sL$ into $M_\sL$ 
determines a presentation of $\pi_1(M_\sL)$ as follows.

Any element of $G_\sL$ can be written as a word
\begin{eqnarray}
\gamma_1\tau_{e_1}\dots\gamma_k\tau_{e_k}\gamma_{k+1}
\end{eqnarray}
satisfying
\begin{description}
\item{(i)} $t(e_i) = i(e_{i+1})$, for $i=1,\dots,k$;
\item{(ii)} $t(e_{k}) = i(e_1) = v$;
\item{(iii)} each $\gamma_i$ is an element of $\pi_1(M_{t(e_i)})$.
\end{description}
The relations on $G_\sL$ are generated by elements of the form
$$
\tau_e\gamma\tau_{\overline{e}} (\gamma')^{-1},
$$
where $\gamma = \phi_e(\kappa)$ and $\gamma' = \phi_{\overline{e}}(\kappa)$
for some $\kappa \in G_e$.
A word of the form (1) is said to have length $k$.   The length
thus depends on the word and not its equivalence class in $G_\sL$.
A word of the form (1) is {\it reduced} if whenever
$e_i = \overline{e_{i+1}}$, for some $i=1,\dots,{k-1}$, then
$\gamma_{i+1} \notin \phi_{e_{i}}(\pi_1(M_{e_i}))$.

The following theorem is useful for determining properties of
groups which are graphs of groups.

\begin{proposition}(\cite{Serre:Trees}, Theorem 11.)
A reduced word of the form (1) 
is trivial if and only if it has length 1 and $\gamma_1 = 1$.
\end{proposition}

For example, we have the following result.

\begin{corollary}  The map
$f : \Gamma_\sL \rightarrow M_\sL$ induces an endomorphism on
fundamental groups.   
\end{corollary}

\heading{Proof.}  Let $G_\sL$ be the image of 
$$
f_* : \pi_1(\Gamma_\sL,v) \rightarrow \pi_1(M_\sL,f(v)).
$$  
Then $G_\sL$ consists of elements of the form
\begin{eqnarray}
f(e_1) \dots f(e_k) = f_* (e_1\dots e_k)
\end{eqnarray}
where $e_1\dots e_k$ is a closed loop on $\Gamma_\sL$
and $i(e_1) = t(e_k) = v$.
Suppose an element of the form (2) is trivial in $G_\sL$
and $k > 1$.   We will show that we can write
the same element in the form (2) with smaller $k$.  
By Theorem 5, $k \ge 2$ and there is at least one $i$ so that 
$e_i = \overline {e_{i+1}}$.  In this case $e_i e_{i+1} = 1$,
so the element can be represented with smaller $k$.  
\qed

We will now present the fundamental groups of $M_\sL$ 
and of $E_\sL$ purely from
the ordered graph $\Gamma_\sL$.

Take any ordered bipartite graph $\widetilde{\Gamma}$ 
with ``line" and ``point" vertices, such that each directed edge
$e(p,L)$ (or $e(L,p)$) joins a point-vertex $v_p$ to a line-vertex
$v_L$ (or vice versa.)  Choose a vertex $v_0$ on $\widetilde{\Gamma}$.
For $v$, a vertex on $\widetilde{\Gamma}$, let
$e_{v,1},\dots,e_{v,d}$ be the ordering of edges emanating from $v$.
Set 
$$
G_v  = \langle\ \lambda_v , \mu_{1,v},\dots,\mu_{d,v} : \lambda_v\mu_{j,v} = 
\mu_{j,v}\lambda_v, \quad j = 1,\dots,d \ \rangle .
$$
(Note that in the previous notation $\lambda_{j,v} = \lambda_v\mu_{j,v}$, 
for each $j = 1,\dots,\degree(v)$.)
Let $\sT_\Gamma$ be the collection of closed loops on $\Gamma$ based at $v_0$.
Then, for any $\tau \in \sT_\Gamma$, $\tau = e_1\dots e_r$, where 
$t(e_j) = i(e_{j+1})$, for $j=1,\dots,r-1$ and
$i(e_1) = t(e_r) = v_0$.  

The following is a corollary of Theorem 2.

\begin{corollary}
Let $\sL$ be a complex line arrangement.
Then $\pi_1(M_\sL)$ can be presented as follows:  
each element of $\pi_1(M_\sL)$ can be written as
$$
\gamma_0 e_1\gamma_1 e_2 \dots e_s \gamma_s,
$$
where $e_1\dots e_s \in \sT_{\Gamma_\sL}$ and $\gamma_j \in G_{t(e_j)}$,
for $j=1,\dots,s$, and the relations are generated by
\begin{eqnarray*}
\alpha e(p,L) \mu_{k,v_L} e(L,p) \beta &=&  
\alpha \lambda_{v_p} \mu_{j,v_p}^{-1}\beta, \\
\alpha e(p,L) \lambda_{v_L} \mu_{k,v_L}^{-1} e(L,p) \beta
&=& \alpha \lambda_{v_p} \beta,
\end{eqnarray*}
if $e(p,L)$ is the $j$th edge emanating from $v_p$,
and $e(L,p)$ is the $k$th edge emanating from $v_L$.
\end{corollary}

Note that the indexing of the edges at each vertex does not
effect the isomorphism class of the group presented.

By Theorem 4, $\pi_1(E_\sL)$ is a quotient of $\pi_1(M_\sL)$
and we have the following presentation of $\pi_1(E_\sL)$. 

\begin{corollary}
Let $\sL$ be a real line arrangement.  Then the fundamental
group of the complement $\pi_1(E_\sL)$ is $\pi_1(M_\sL)$
with the additional relations
$$
f(e_1)\dots f(e_s) = 1
$$
where $f$ is as in Theorem 4 and $e_1 \dots e_s \in \sT_\Gamma$. 
\end{corollary}

\section{Complex line arrangements}

Let $\sL$ be an arbitrary complex line arrangement consisting
of lines $\sA$ and points of intersection $\sP$.  
The proofs of Theorem 3 and Theorem 4 given in Section 3 rely
on the existence of a continuous map
$$
g : \Sigma_\sL \rightarrow M_\sL
$$
from the skeleton $\Sigma_\sL$ of $\sL$ to the boundary 3-manifold $M_\sL$
whose image contracts in $E_\sL$.

We will show that Theorem 3 and Theorem 4 cannot be applied
to arbitrary complex line arrangements.  To do this we 
extend the definition of skeleton for complex
line arrangements, and show that the map  $g$ and hence the map
$$
f : \Gamma_\sL \rightarrow M_\sL
$$
with the required properties does not exist in general.

\heading{Knotted 1-complexes.}

Let $\Sigma \subset \R^3$ be a compact 1-complex.  We will
call $\Sigma$ a {\it graph knot}.  
Let $N_\Sigma$ be a thickening of $\Sigma$ in $\R^3$  and let
$M_\Sigma$ be its boundary.  Then $N_\Sigma$ is a
handlebody and $M_\Sigma$ is an oriented surface.
Let $\alpha' : N_\Sigma \rightarrow \Sigma$ be the contraction
map and let $\alpha: M_\Sigma \rightarrow \Sigma$ be the 
restriction of $\alpha'$ to $M_\Sigma$.
We say $\Sigma$ is {\it unknotted} if 
there is an embedding
$$
g : \Sigma \rightarrow M_\Sigma
$$
such that $\alpha\circ g$ is the identity on $\Sigma$,
and the image of $g$ contracts in $\R^3$.
We say that $\Sigma$ is {\it knotted} otherwise.  

Dehn's Lemma (\cite{Rolfsen:Knots}, p. 101) implies the following.

\begin{lemma}  A graph knot $\Sigma \subset \R^3$ is
knotted if there is an embedding
$$
\sqcup S^1 \rightarrow \Sigma
$$
whose image is a nontrivial knot or link in $\R^3$.
\end{lemma}

\heading{Skeleta for complex arrangements.}

Let $\sL$ be an arbitrary complex line arrangement and 
let
$$
\rho : \C^2 \rightarrow \C
$$
be projection onto the $x$-coordinate.  Assume 
that $\rho|_L$ is one-to-one for all lines $L \in \sA$.
Let $\sQ = \rho(\sP)$.

Given an ordering of the points $q_1\dots,q_s \in \sQ$,
let $I = I_1 \cup \dots \cup I_{s-1}$ be a 
union of 
embedded line segments in $\C$ such that each $I_i$ has endpoints
$q_i$ and $q_{i-1}$ and no pair $I_i$ and $I_j$ intersect 
except possibly at their endpoints.  We will call
$\Sigma_\sL(I) = \rho^{-1}(I) \cap \sL$ the {\it skeleton}
of $\sL$ associated to $I$. 
(The 1-complex $\Sigma_\sL(I)$  considered as a subset of
$\R^2 \times I$ is also known as a {\it braided wiring
diagram} and has been used to
describe the fundamental group and homotopy type of the complement
of arbitrary complex line arrangements and algebraic plane curves
\cite{Arv:Fund}, \cite{C-S:Braided}.)

The skeleton
$\Sigma_\sL$ defined for real line arrangements comes from ordering
the $q_i$ in the order in which they lie on the real line $\R$ and
taking $I$ to be the line segment on $\R$ connecting the smallest
to the largest.
We will call this $\Sigma_\sL$ the {\it standard skeleton} for the
real line arrangement $\sL$.

As with real line arrangements, we have the following.

\begin{lemma} For any complex line arrangement $\sL$ and any ordering
of $\sQ$, the 1-complex $\Sigma_\sL(I)$
is homotopy equivalent to the graph $\Gamma_\sL$.   
\end{lemma}

We can think of $\Sigma_\sL(I)$ as a subset of $\R^3 = \R^2 \times \R$
by embedding $I \in \R$.  This makes  $\Sigma_\sL(I)$ a graph knot. 
When $\sL$ is a real line arrangement and $\Sigma_\sL$
is the standard skeleton $\Sigma_\sL$ is unknotted.
Conversely, if $\sL$ is a complex line arrangement such that for
all orderings the associated skeleton $\Sigma_\sL(I)$ is knotted,
then $\sL$ is not topologically equivalent to a real line arrangement.

Lemma 3.1 then generalizes as follows.

\begin{theorem} Let $\sL$ be a complex line arrangement
and let $\alpha : M_\sL \rightarrow \sL$ be the natural
projection map of the boundary manifold $M_\sL$ onto $\sL$.
Then for some choice of $I$, $\Sigma_\sL(I)$ is unknotted
if and only if there is an embedding 
$$
g : \Sigma_\sL(I) \rightarrow M_\sL,
$$
such that   
\begin{description}
\item{(i)} $\alpha \circ g$ is the identity map; and
\item{(ii)} the complement $E_\sL$ is homotopy 
equivalent to 
$$
M_\sL /g(\Sigma_\sL(I)).
$$
\end{description}
\end{theorem}

\heading{Proof.}
The proof is the same as for Lemma 3.1.
\qed

It follows that if $\Sigma_\sL(I)$ is knotted for all choices of $I$,
then Theorem 4 cannot be extended to $\sL$.  For general algebraic
plane curves, one can also describe the boundary 3-manifold $M_\sC$ as   
a graph of manifolds over a suitably defined incidence graph 
$\Gamma_\sC$, but Theorem 5 provides an obstruction for mimicing 
Theorem 4 in this setting.

\end{document}